\documentclass[12pt]{amsart}
\usepackage{amssymb,amscd,amsmath}
\usepackage{amsfonts}
\usepackage{fullpage}
\usepackage{array}
\usepackage{amsthm}
\usepackage{url}
\usepackage[top=1.18in, left=1in, right=1in, bottom=1.8in]{geometry}
\newtheorem{theorem}{Theorem}[section]

\newtheorem{lemma}[theorem]{Lemma}

\newtheorem{conjecture}[theorem]{Conjecture}
\theoremstyle{definition}
\newtheorem{definition}[theorem]{Definition}
\newtheorem{remark}[theorem]{Remark}

\def\GL{\operatorname{GL}}

\input xy
\xyoption{all}
\title{Notes on Low Degree $L$-Data}
\author{Thomas Oliver}
\thanks{Heilbronn Institute for Mathematical Research, University of Bristol, Bristol, UK. \\ email: to14172@bristol.ac.uk}
\date{\today}
\begin{document}
\maketitle
\subsection*{Abstract} These notes are an extended version of a talk given by the author at the conference ``Analytic Number Theory and Related Areas'', held at Research Institute for Mathematical Sciences, Kyoto University in November 2015. We are interested in ``$L$-data'', an axiomatic framework for $L$-functions introduced by Andrew Booker in 2013 \cite{LFAD}. Associated to each $L$-datum, one has a real number invariant known as the degree. Conjecturally the degree $d$ is an integer, and if $d\in\mathbb{N}$ then the $L$-datum is that of a $\GL_n(\mathbb{A}_F)$-automorphic representation for $n\in\mathbb{N}$ and a number field $F$ (if $F=\mathbb{Q}$, then $n=d$.). This statement was shown to be true for $0\leq d<\frac{5}{3}$ by Booker in his pioneering paper \cite{LFAD}, and in these notes we consider an extension of his methods to $0\leq d<2$. This is simultaneously a generalisation of Booker's result and the results and techniques of Kaczorowski--Perelli in the Selberg class \cite{OTSOSCVII}. Furthermore, we consider applications to zeros of automorphic $L$-functions. In these notes we review Booker's results and announce new ones to appear elsewhere shortly \cite{COLDLD}.
\subsection*{Acknowledgement}
The author is very happy to have had the opportunity to give this lecture and is grateful for the invitation and help of the organiser of the conference, Professor Yuichi Kamiya. Moreover, the author appreciates the efforts and assistance of Masatoshi Suzuki regarding the organisation of his time in Japan - surely things would have gone much less smoothly without him! Finally, sincere thanks are made to Andrew Booker, who first suggested this direction of research. The author was supported by a Heilbronn postdoctoral research fellowship.
\section{Introduction}
The \textsl{Selberg class} is an axiomatic framework for $L$-functions, introduced by Selberg in 1989 \cite{OANCARAACODS}. Specifically, it is the set of complex functions $L(s)$ satisfying the following 5 axioms:
\begin{enumerate}
\item \textbf{Dirichlet series} - there are $a_n\in\mathbb{C}$ such that $L(s)=\sum_{n=1}^{\infty}a_nn^{-s}$, with absolute convergence for $\Re(s)>1$;
\item \textbf{Analytic continuation} - there is $m\in\mathbb{Z}_{\geq0}$ such that $(s-1)^mL(s)$ continues to an entire function of finite order;
\item \textbf{Functional equation} - there are $k\in\mathbb{Z}_{\geq0}$, $Q\in\mathbb{R}_{>0}$, $\lambda_1,\dots,\lambda_k\in\mathbb{R}_{>0}$, $\mu_1,\dots,\mu_k\in\{\Re(z)>0\}$, $\epsilon\in\{z\in\mathbb{C}:|z|=1\}$ such that
\[\Lambda(s)=\overline{\Lambda(1-\bar{s})},\text{ where }\Lambda(s)=\epsilon Q^sL(s)\prod_{j=1}^k\Gamma(\lambda_js+\mu_j);\]
\item \textbf{Ramanujan hypothesis} - For every $\varepsilon>0$, $a_n\ll_{\varepsilon}n^{\varepsilon}$;
\item \textbf{Euler product} - $a_1=1$ and $\log L(s)=\sum_{n=2}^{\infty}b_nn^{-s}$, where $b_n$ is supported on prime powers and $b_n\ll n^{\theta}$ for some $\theta<\frac{1}{2}$.
\end{enumerate}
This set of axioms captures the sort of behaviour believed to be necessary for $L(s)$ to satisfy the (generalized) Riemann hypothesis. On the other hand, each Dirichlet series in the Selberg class is supposed to be an \textsl{automorphic} $L$-function. In turn, such $L$-functions  constitute a rather large sample of all known $L$-functions.  Indeed, many familiar objects such as Dirichlet characters, modular forms, Maass forms etc. give rise to automorphic representations. Additionally, the so-called \textsl{Langlands philosophy} would have it that all \textsl{motivic} $L$-functions arising from arithmetic geometry are automorphic. For example, this includes the Artin $L$-functions associated to complex representations of Galois groups, and Hasse--Weil $L$-functions associated to $l$-adic representations of Galois groups constructed from algebraic varieties over number fields. In general, automorphic $L$-functions are not known to satisfy the Ramanujan hypothesis\footnote{Contrast this to the \textsl{Selberg orthogonality conjecture}, which is close to being settled for automorphic $L$-functions. It is worth noting that automorphic $L$-functions are in the \textsl{extended} Selberg class, which is defined to be those functions satisfying only the analytic axioms (1) - (3). Kaczorowski and Perelli have managed to classify low degree elements of not only the Selberg class, but moreover low degree elements in this extension \cite{OTSOSCVII}.}.
\newline Motivated by both philosophical and practical considerations, Booker introduced an axiomatic framework for the study of automorphic $L$-functions \cite{LFAD}. Booker's basic idea was to parametrise \textsl{explicit formulae}, of the type introduced by Weil \cite{SLFEDLTDNP}.  His class, referred to here as the class of $L$-\textsl{data}, includes not only the Selberg class, but also the class of automorphic $L$-functions, and even, the Artin $L$-functions. An advantage of the class of $L$-data over the Selberg class is the immediate applications to the study of vanishing orders of automorphic $L$-functions and the cancellation of zeros between different automorphic $L$-functions. This is down to a flexibility in the admissible gamma factors which is not present in Selberg's axioms (cf. section 2).
\newline Theorems concerning zeros of automorphic $L$-functions follow immediately from the classification of \textsl{positive} $L$-data. The term ``positive'' can be interpreted loosely as having finitely many poles. Extending the theory for the Selberg class, this classification is built on an invariant called the \textsl{degree}. A priori the degree is a real number, though conjecturally it is in fact integral. Moreover, if the degree of an $L$-datum is an integer, one expects that the $L$-datum corresponds to the $L$-function of a $\GL_n(\mathbb{A}_F)$-automorphic representation for some number field $F$ and some $n\in\mathbb{N}$. Statements of this nature are often referred to as \textsl{converse theorems}. If $F=\mathbb{Q}$, then $n$ is in fact the degree of the $L$-datum. For example, each finite order Hecke charcter of $\mathbb{Q}$ (that is $\GL_1(\mathbb{A}_{\mathbb{Q}})$-automorphic representation) corresponds to a unique primitive Dirichlet character $\chi$, and it is known that every degree 1 $L$-datum arises from the $L$-function of such a $\chi$. On the other hand, a degree 2 $L$-function could arise from not only a $\GL_2(\mathbb{A}_{\mathbb{Q}})$-automorphic representation, such as a modular form or a Maass form, but also a $\GL_1(\mathbb{A}_F)$-automorphic representation (Hecke character) for some quadratic extension $F$ of $\mathbb{Q}$ (in the case of the trivial Hecke character we have a product of two Dirichlet $L$-functions). In these notes we will work only over $\mathbb{Q}$. 
\newline From the classification of positive $L$-data of degree $0\leq d<2$ one can uniformly prove a wide family of theorems concerning the vanishing order of automorphic $L$-functions. To show the malleability of the method, we state a whimsical example of such a result: If $\pi$ is a unitary cuspidal  $\GL_{163}(\mathbb{A}_{\mathbb{Q}})$-automorphic representation, then the completed $L$-function $\Lambda(s,\pi)$ has infinitely many zeros of order not divisible by $82$. Indeed, once the classification of $L$-data is settled, all one has to check to prove this is that $163/82<2$. That statement is rather convoluted, but the very same logic dictates that the completed $L$-function $\Lambda(s,\rho)$ of a unitary cuspidal automorphic representation $\rho$ of $\GL_3(\mathbb{A}_{\mathbb{Q}})$ has infinitely many zeros of odd order. This time, the key is that $3/2<2$. 
\newline The contents of these notes break down as follows. Section 2 is devoted to the formal definitions of the concepts discussed above. Section 3 then discusses some theorems (old and new) regarding the classification of low degree $L$-data. Finally section 4 is concerned with sketching the proofs, with precise details left to the references given there.
\section{Definitions}
First we remind the reader of explicit formulae - these are distribution identities which relate zeros of $L$-functions to sums over primes. Of particular importance to us will be automorphic $L$-functions, so we focus on this case. Let $\pi$ be a unitary cuspidal automorphic representation of $\GL_d(\mathbb{A}_{\mathbb{Q}})$, with $L$-function $L(s,\pi)$ and conductor $q$. There are numbers $\mu_j\in\mathbb{C}$ such that the \textsl{completion} $\Lambda(s)$ is defined as follows
\[\Lambda(s,\pi):=L(s,\pi_{\infty})L(s,\pi),\]
where 
\[L(s,\pi_{\infty}):=\prod_{j=1}^d\Gamma_{\mathbb{R}}(s+\mu_j);~\Gamma_{\mathbb{R}}(s)=\pi^{-s/2}\Gamma(s/2).\]
We may write the logarithmic derivative of $L(s,\pi)$ as a Dirichlet series:
\[-\frac{L'}{L}(s,\pi)=\sum_{n=2}^{\infty}c_nn^{-s}.\]
If $g:\mathbb{R}\rightarrow\mathbb{C}$ is a smooth function of compact support with Fourier transform $h(z)=\int_0^{\infty}g(x)e^{ixz}dx$ such that $h(\mathbb{R})\subseteq\mathbb{R}$, then
\[\sum_{z\in\mathbb{C}}\text{ord}_{s=\frac{1}{2}+iz}\Lambda(s,\pi)\cdot h(z)\]
\[=2\Re[\int_0^{\infty}(g(0)-g(x))\sum_{j=1}^d\frac{e^{-(\frac{1}{2}+\mu_j)x}}{1-e^{-2x}}dx-g(0)(\frac{1}{2}\log q-\Re\sum_{j=1}^d\frac{\Gamma_{\mathbb{R}}'}{\Gamma_{\mathbb{R}}}(\frac{1}{2}+\mu_j))-\sum_{n=2}^{\infty}\frac{c_ng(\log(n))}{\sqrt{n}}]\]
Note that such formulae are \textsl{additive} in the sense that the formula for a product of $L$-functions is a sum of the formulae for the factors. Later we will see that this allows us to study vanishing orders of automorphic $L$-functions via the most basic linear algebra. The additivity makes it clear that if we can encapsulate the essence of explicit formulae, we can, for example, incorporate quotients of automorphic $L$-functions into our framework (as a difference of explicit formulae). A related point of view is that this additivity allows us to deform the gamma factors and put them on the same footing as the non-Archimedean Euler factors. We do this via the integral kernel, for example, halving it gives us the square root of the gamma factor. Neither of these features appear in the Selberg class - the functional equation of a square root or quotient of $L$-functions does not have gamma factors of the correct form. On the other hand, according to the Euler product axiom, a non-Archimedean Euler factor in the Selberg class can be any function of the form $e^{f(p^{-s})}$ where $f$ is analytic on a disc of radius $p^{-\theta}$, $\theta<\frac{1}{2}$. We remark that, whilst in both the archimedean and non-archimdean case the Euler factors are more flexible than what seems to occur in nature, it is the lack of uniformity than concerns us most. 
\newline Bearing this and more in mind, Booker suggested the definition below in 2013 \cite{LFAD}.
\begin{definition}\label{Ldata.definition}
An $L$-datum is a triple $F=(f,K,m)$, where
\[f:\mathbb{Z}_{>0}\rightarrow\mathbb{C};~K:\mathbb{R}_{>0}\rightarrow\mathbb{C};~m:\mathbb{C}\rightarrow\mathbb{R};\]
are such that
\begin{enumerate}
\item\textbf{Growth} - $f(1)\in\mathbb{R}$, $f(n)\log^kn\ll_k1$ for all $k>0$, and $\sum_{n\leq x}|f(n)|^2\ll_{\varepsilon}x^{\varepsilon}$, for all $\varepsilon>0$;
\item\textbf{Degree} - $xK(x)$ extends to a Schwartz function on $\mathbb{R}$ and $\lim_{x\rightarrow0^{+}}xK(x)\in\mathbb{R}$;
\item\textbf{Multiplicity} - supp$(m):=\{z\in\mathbb{C}:m(z)\neq0\}$ is discrete and contained in a horizontal strip $\{z\in\mathbb{C}:|\Im(z)|\leq y\}$ for some $y\geq 0$. Moreover $\sum_{z\in\text{supp}(m),|\Re(z)|\leq T}|m(z)|\ll1+T^A$ for some $A\geq0$ and $\#\{z\in\text{supp}(m):m(z)\notin\mathbb{Z}\}<\infty$;
\item\textbf{Explicit formula} - For every smooth function $g:\mathbb{R}\rightarrow\mathbb{C}$ of compact support and Fourier transform $h(z)$ satisfying $h(\mathbb{R})\subseteq\mathbb{R}$ we have the equality
\[\sum_{z\in\text{supp}(m)}m(z)h(z)=2\Re[\int_0^{\infty}K(x)(g(0)-g(x))dx-\sum_{n=1}^{\infty}f(n)g(\log n)].\]
\end{enumerate}
Given an $L$-datum $F=(f,K,m)$, we define the \textsl{L-function} of $F$ to be the following Dirichlet series.
\[L_F(s):=\exp(\sum_{n=2}^{\infty}\frac{f(n)}{\log n}n^{\frac{1}{2}-s})=:\sum_{n=1}^{\infty}a_F(n)n^{-s},~\Re(s)>1.\] 
The \textsl{degree} $d_F$ of $F$ is defined to be
\[d_F:=2\lim_{x\rightarrow0^{+}}xK(x).\] 
The \textsl{analytic conductor} $Q_F$ of $F$ is defined to be
\[Q_F=e^{-2f(1)}.\]
We say that $f$ is \textsl{positive} if there are at most finitely many $z\in\mathbb{C}$ with $m(z)<0$.
\end{definition}
\begin{remark}
On first reading, it might seem somewhat strange to allow for the possibility of finitely many $z\in\mathbb{C}$ such that $m(z)\notin\mathbb{Z}$. Afterall, $m$ is supposed to act like the order of a zero or pole at $z$. What pushes us in this direction is that in applications one sometimes seeks to scale $L$-data by fractional constants. This is very useful in the study of vanishing orders of automorphic $L$-functions.
\end{remark}
We have seen at the start of this section that the unitary cuspidal automorphic representations of $\GL_n(\mathbb{A}_\mathbb{Q})$ give rise to \textbf{positive} $L$-data. Moreover, the degree of the $L$-datum of a unitary cuspidal $\GL_n(\mathbb{A}_{\mathbb{Q}})$-automorphic representation is easily calculated to be $n$. Specifically, if $\pi$ has conductor $q$ then, in the notation from the start of this section, one has
\[F_{\pi}:=(f_{\pi},K_{\pi},m_{\pi});\]
where
\[f_{\pi}(n):=\begin{cases}
-\frac{1}{2}\log q-\Re\sum_{j=1}^d\frac{\Gamma'_{\mathbb{R}}}{\Gamma_{\mathbb{R}}}(\frac{1}{2}+\mu_j), & n=1;\\
\frac{c_n}{\sqrt{n}}, & n>1;
\end{cases}\]
\[K_{\pi}(x):=\sum_{j=1}^d\frac{e^{-(\frac{1}{2}+\mu_j)x}}{1-e^{-2x}};\]
\[m_{\pi}(z):=\text{ord}_{s=\frac{1}{2}+iz}\Lambda(s,\pi).\]
That these functions satisfy all the axioms in definition~\ref{Ldata.definition} is explained in \cite[Example~1.4]{LFAD}. Similarly, Artin $L$-functions define $L$-data and the conjecture that Artin $L$-functions are automorphic amounts to the statement that the associated $L$-data are positive. One may also show that Dirichlet series in the Selberg class give rise to $L$-data. It could be helpful to keep in mind the diagram of (conjectural) inclusions below.
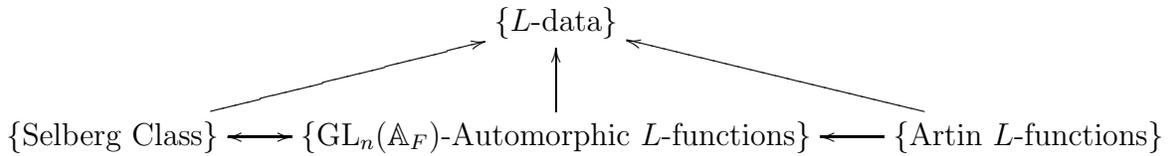
\[\xymatrix{ 
 & \{L\text{-data}\} & \\
\{\text{Selberg Class}\}\ar@{->}[ur]\ar@{->}[r] & \{\GL_n(\mathbb{A}_F)\text{-Automorphic }L\text{-functions}\}\ar@{->}[l]\ar@{->}[u] & \{\text{Artin }L\text{-functions}\}\ar@{->}[ul]\ar@{->}[l] \\
}\]
We offer the following caveats:
\begin{enumerate}
\item It is not yet clear whether automorphic $L$-functions are in the Selberg class and vice versa, though there are conjectures in this direction as discussed.
\item Artin $L$-functions are known to admit \textsl{meromorphic} continuation to $\mathbb{C}$ but are not yet known to be automorphic in general. For a discussion of this see \cite[chapter~4]{ITLP}.
\item Hasse--Weil $L$-functions do not seem to fit neatly into this picture yet. In very special cases there are automorphicity results (cf. [loc. cit., chapter~5]).
\end{enumerate}
Positive $L$-data should be classified by their degree. Specifically we formulate the following:
\begin{conjecture}\label{classification.conjecture}
If $F$ is a positive $L$-datum, then the degree $d_F\in\mathbb{N}$. Moreover, the $L$-function $L_F(s)$ of $F$ is (up to an imaginary displacement) the $L$-function of a $\GL_n(\mathbb{A}_F)$-automorphic representation for some number field $F$ and $n\in\mathbb{N}$.
\end{conjecture}
This is an extension of the analogous conjecture in the Selberg class. The current best result in that setting is a complete classification for $0\leq d<2$ which was published by Kaczorowski--Perelli in 2011 \cite{OTSOSCVII}. Indeed, 2 is something of a natural boundary as in this degree one encounters a whole host of as yet mysterious $L$-functions, eg. those of elliptic curves; modular forms; and Maass forms. The methods of Kaczorowski--Perelli offer a complexity-type argument for this heuristic.
\newline In degree 0, we have only the \textsl{trivial} $L$-data. Formally, this is a ``multiplicity one'' statement.
\begin{theorem}\textbf{\cite[theorem~1.6]{LFAD}.}
For an $L$-datum $F=(f,K,m)$, we have
\[F=(0,0,0)\Leftrightarrow\sum_{n=2}^{\infty}\frac{|f(n)|}{\log n}<\infty\Leftrightarrow\sum_{n=1}^{\infty}\frac{|a_F(n)|}{\sqrt{n}}<\infty\]
\[\Leftrightarrow L_F(s)\text{ is a ratio of Dirichlet polynomials}\Leftrightarrow\sum_{z\in\text{supp}(m),|\Re(z)|\leq T}|m(z)|=o(T).\]
\end{theorem}
In particular, if $d_F=0$, then $L_F(s)=1$ as expected. The theorem above teaches us that we can think about $L$-data as Dirichlet series without loss of crucial information. These $L$-functions are shown to have analogues of the usual analytic properties in \cite[Proposition~2.1]{LFAD} (cf. section 4). 
\newline Multiplicity 1 can be used moreover to deduce the non-existence of $L$-data of degree $0<d<1$ - all that is required is Mellin inversion and Stirling's formula \cite[section~3.1]{LFAD}. Furthermore, we can say that if $F$ is a positive $L$-data of degree 1, then there is a Dirichlet character $\chi$ and $t\in\mathbb{R}$ such that 
\[L_F(s)=L(s+it,\chi).\]
As Dirichlet characters are precisely the finite order Hecke characters of $\mathbb{Q}$, ie. $\GL_1(\mathbb{A}_{\mathbb{Q}})$-automorphic representations, this is consistent with conjecture~\ref{classification.conjecture}. The proof, explained in \cite[section~3.2]{LFAD}, works by firstly showing that the coefficients of $L_F(s)$ are periodic and then applying a result of Saias--Weingartner \cite{ZODSWPC}, which gives conditions under which Dirichlet series with periodic coefficients arise from Dirichlet characters. The periodicity of the coefficients follows from the reflection formula for the gamma function.
\newline As an application of the classification of degree $1$ $L$-data, one can answer particular questions about cancellation of zeros between automorphic $L$-functions. For example, let $\pi_1$ and $\pi_2$ be non-isomorphic unitary cuspidal automorphic representations for $\GL_{d_1}(\mathbb{A}_{\mathbb{Q}})$ and $\GL_{d_2}(\mathbb{A}_{\mathbb{Q}})$ respectively. If $d_2-d_1\leq1$, then the quotient $\Lambda(s,\pi_2)/\Lambda(s,\pi_1)$ has infinitely many poles. Indeed, if the quotient has only finitely many poles, then the associated $L$-data is positive and we have either the trivial $L$-data or the $L$-data of a Dirichlet $L$-function. The first option violates the \textsl{non-isomorphic} condition and the second violates the \textsl{cuspidal} condition. There is precedent for this line of enquiry - similar results have been proved by Raghunathan in special cases \cite{ACOZOLF}, \cite{CTDSP}, \cite{OLFWPSMFE}. For example, it was shown in \cite{ACOZOLF} that the following quotient has infinitely many poles
\[\frac{L(\text{Sym}^2(\pi_f)\otimes\chi,s)}{L(\chi,s)},\]
where $\pi_f$ is the cuspidal automorphic representation associated to a cuspidal modular form $f$. Note that here the difference in degree is 2. More generally, questions concerning cancellation of zeros can be couched in terms of the \textsl{Grand Simplicity Hypothesis}, which concerns linear independence of zeros. In the next section we will discuss further conjectures concerning cancellation of zeros, especially when the degrees differ by 2.
\begin{remark} 
One may relax the cuspidality assumption. This allows the statement to be formulated for products of cuspidal $L$-functions. Quotients of products of automrophic $L$-functions arise (at least conjecturally) as the zeta functions of arithmetic schemes.
\end{remark}
\section{Degrees $1<d\leq2$}
Extending the classification of positive $L$-data to degrees $d>1$ will allow us to deduce more facts about zeros of $L$-functions. For example, if we knew that there were no $L$-data of degree $1<d<\frac{3}{2}+\varepsilon$, $\varepsilon>0$, then it would follow that the completed $L$-function $\Lambda(s,\pi)$ of a unitary cuspidal automorphic representation $\pi$ of $\GL_3(\mathbb{A}_{\mathbb{Q}})$ has infinitely many zeros of odd order. Indeed, let the positive $L$-datum $F$ be that associated to $\pi$. If $\Lambda(s,\pi)$ has at most finitely many zeros of odd order, then $m(z)$ is an even integer for almost all $z$ and $\frac{1}{2}F$ is a positive $L$-datum of degree $\frac{3}{2}$. There are a few observations to make about this argument:
\begin{itemize}
\item It is clear that similar results for higher degree automorphic representations would follow from suitable non-existence results, though beyond degree $2$ is out of reach; 
\item One can use the same argument to prove statements such as \textsl{the }$L$\textsl{-function of a cuspidal automorphic representation of degree 4 has infinitely many zeros of order not divisible by 3}. The statement about zeros of \textsl{odd} order is of greater historical significance, especially results for simple zeros.
\end{itemize}
In \cite[theorem~1.7]{LFAD}, Booker proved the non-existence result for $1<d<\frac{5}{3}$. This is sufficient for the $\GL_3(\mathbb{A}_{\mathbb{Q}})$ argument. The author has subsequently extended the techniques of Kaczorowski--Perelli \cite{OTSOSCVII} to $L$-data to obtain:
\begin{theorem}\label{d<2.theorem}
There are no positive $L$-data of degree $1<d<2$.
\end{theorem}
The proof of this is to appear in \cite{COLDLD} (cf. section 4 for a limited sketch). Jointly with Michael Neururer, the author is taking tentative steps with degree $2$ $L$-data. In this setting, one would like a converse theorem of the form ``if a Dirichlet series has degree 2 and \textbf{nice} analytic properties, then it is the $L$-function of a modular form; the $L$-function of a Maass form; or a Hecke $L$-function for a quadratic number field''. Moreover, each category should be characterised by the associated functional equation. The word \textbf{nice} should allow for poles but need not necessarily involve an Euler product. Though there are relevant representation-theoretic converse theorems, there is not yet such a general converse theorem for degree 2 Dirichlet series which does not assume an Euler product\footnote{There are a few noteworthy remarks concerning the literature on this. The first thing to acknowledge is that there is an extensive literature on representation-theoretic converse theorems relevant to the case of degree 2 Dirichlet series written in the adelic language, for example \cite{CTFGLn}, \cite{CTFGLnII}. In this setting, the Euler product of automorphic $L$-functions is self-evident. On the other hand, Weil's original converse theorem for Dirichlet series coming from homolorphic modular forms made no use of the Euler product \cite{UDBDRDF}. Booker and Krishnamurthy have proved a generalisation of this allowing for the Dirichlet series to have a wider class of poles \cite{WCTWP} (see also \cite{CTDSP}). A Weil-type converse theorem for the Dirichlet series of Maass forms was stated in \cite{BockleThesis}, though there is an apparently undocumented error in the statement of the non-holomorphic analogue of the fact that if a holomorphic function on the upper half-plane is invariant under an elliptic operator then it is constantly zero (cf. \cite[Lemma~1.5.1]{AFAR}). More recently, Raghunathan has proved a converse theorem for Dirichlet series with certain poles satisfying Maass's functional equation \cite{OLFWPSMFE}.}. On the other hand, with a converse theorem of this nature in mind, one conjectures statements along the following lines. 
\begin{conjecture}
Let $\pi_1$, $\pi_2$ be non-isomorphic unitary cuspidal automorphic representations of $\GL_{d_1}(\mathbb{A}_{\mathbb{Q}})$ and $\GL_{d_2}(\mathbb{A}_{\mathbb{Q}})$, respectively. If $d_2-d_1\leq 2$, then the quotient $\Lambda(s,\pi_2)/\Lambda(s,\pi_1)$ has infinitely many poles.
\end{conjecture}
The proof of this would follow much as in the case $d_2-d_1=1$ presented in the previous section. Note that Ragunathan's results found in  \cite{ACOZOLF}, \cite{CTDSP}, \cite{OLFWPSMFE} are all deduced from appropriate degree 2 converse theorems. Cases of the above conjecture will be proved elsewhere in collaboration with Neururer. 
\section{Towards A Proof}
We conclude these notes with a brief outline of the proof of theorem~\ref{d<2.theorem}. The theory is somewhat intricate and constraints on the length of contributions to these proceedings mean that we cannot go into great detail. A full proof will appear in \cite{COLDLD}.
\newline Assume for a contradiction that $F$ is an $L$-datum of degree $1<d<2$ with associated $L$-function $L_F(s)=\sum_{n=1}^{\infty}a_F(n)n^{-s}$. Inspired by Booker \cite{LFAD} (see also references therein), the proof of theorem~\ref{d<2.theorem} is based on variants of the following natural exponential sum:
\[S_F(z):=\sum_{n=1}^{\infty}a_F(n)\exp(2\pi inz).\]
Studying the behaviour of this sum leads to the result for $1<d<\frac{5}{3}$ as in \cite{LFAD}. To push this technique to $d<2$, inspired by the proof of the analogous result in the Selberg class by Kaczorowski--Perelli \cite{OTSOSCVII}, we consider the more general sum
\[S_F(z;\underline{\alpha})=\sum_{n=1}^{\infty}a_F(n)\exp(2\pi i(c_1n^{\alpha_1}+\dots+c_Nn^{\alpha_N})z);\]
\[\underline{\alpha}:=(\alpha_1,\dots,\alpha_N),~(c_1,\dots,c_n)\in(\mathbb{R}_{>0})^N.\]
Of particular interest is the case $N=1$, $\alpha_1=\frac{1}{d}$. The expression $S_F(z;\underline{\alpha})$ above should be compared to to so-called \textsl{multidimensional non-linear twists} of Kaczorowski--Perelli \cite{OTSOSCV}, \cite{OTSOSCVI}, \cite{OTSOSCVII}. To proceed, one needs to understand the analytic properties of the Dirichlet series $L_F(s)$ associated to the $L$-datum $F$. In fact, by reversing the steps taken in the proof of the explicit formula one can show that the $L$-function of an $L$-datum $F=(f,K,m)$ admits meromorphic continuation and functional equation much like Dirichlet series in the Selberg class. Specifically, there is a function $\gamma_F(s)$ defined uniquely up to real scalars such that
\begin{itemize}
\item $\log\gamma_F(s)$ is holomorphic for $\Re(s)>\frac{1}{2}$, and $\frac{d^n}{ds^n}\log\gamma_F(s)$ extends continuously to $\Re(s)\geq\frac{1}{2}$ for each $n\geq0$;
\item there are constants $d,c_{-1}\in\mathbb{R}$ and $\mu,c_0,c_1,\dots\in\mathbb{C}$ such that
\[\log\gamma_F(s)=(s-\frac{1}{2})(\frac{d}{2}\log\frac{s}{e}+c_{-1})+\frac{\mu}{2}\log\frac{s}{2}+\sum_{j=0}^{n-1}\frac{c_j}{s^j}+O_n(|s|^{-n}),\]
uniformly for $\Re(s)\geq\frac{1}{2}$ and any fixed $n\geq0$. In fact, the number $d$ turns out to be the degree;
\item The product $\Lambda_F(s)=\gamma_F(s)L_F(s)$ continues meromorphically to 
\[\Omega=\mathbb{C}-\{\text{finitely many vertical rays}\}\]
\[-\{\text{finitely many horizontal line segments}\}.\]
Moreover $\Lambda_F(s)$ has meromorphic finite order on $\Omega$, that is the ratio of two finite order holomorphic functions on $\Omega$. The meromorphy fails to extend to $\mathbb{C}$ on account of the fact that $m$, supposed to play the role of a multiplicity, need not be integral at all points;
\item The functional equation $\Lambda_F(s)=\overline{\Lambda_F(1-\bar{s})}$ holds as an identity of meromorphic functions on $\Omega$;
\item The logarithmic derivative of $\Lambda_F(s)$ continues meromorphically to $\mathbb{C}$, with at most simple poles, and satisfies
\[\text{Res}_{s=\frac{1}{2}+iz}\frac{\Lambda_F'}{\Lambda_F}(s)=m(z).\]
In particular supp$(m)\subseteq\{z\in\mathbb{C}:|\Im(z)|\leq\frac{1}{2}\}$.
\end{itemize}
Using Stirling's formula and Mellin inversion, one can build functions $G(s;\alpha_i)$ from the gamma function such that the $k$-th derivative of $S(z;\underline{\alpha})$ is a sum of iterated inverse Mellin transforms. Specifically:
\[z^kS_F^{(k)}(z;\underline{\alpha})=\]
\[O(\Im(z)^{-\varepsilon})+\sum_{j=0}^k\frac{c_{kj}}{(2\pi i)^{N+1}}\int_{\Re(s_1)=\sigma_1}\dots\int_{\Re(s_N)=\sigma_N}\Lambda_F(\sum_{i=0}^Ns_i)\prod_{i=0}^NG(s_i+\delta(\alpha_i,j);\alpha_i)(-iz)^{-\frac{s_i}{\alpha_i}}ds_i;\]
for some constants $c_{kj}\in\mathbb{C}$, $0\leq j\leq k$, $c_{kk}\neq0$, where
\[\delta(\alpha_i,k):=\frac{2j\alpha_i}{2-d\alpha_i},\]
and the contours are chosen in order that the integral be well-defined. For example, if $N=1$ and $\alpha_1=1$, then
\[G(s)=(2\pi(1-\frac{d}{2})^{1-\frac{d}{2}}e^{c_{-1}})^{\frac{1}{2}-s}\Gamma((1-\frac{d}{2})(s-\frac{1}{2})+\frac{1-\mu}{2}),\]
where the constants $c_{-1}\in\mathbb{R}$ and $\mu\in\mathbb{C}$ come from the analytic properties of $L_F(s)$. 
In general, the iterated contour integral is hard to deal with. That said, we can reduce the ``depth'' by the residue theorem. Let $[N]=\{0,\dots,N\}$, then
\[S_F(z;(\alpha_0,\dots,\alpha_N))=O_{k,\varepsilon}(\Im(z)^{-\varepsilon})+L_F(N)\]
\[+\sum_{\emptyset\neq\mathcal{A}\subset[N]}\frac{1}{(2\pi i)^{|\mathcal{A}|}}\int_{\Re(s_{n_1})=\tau_{n_1}}\cdots\int_{\Re(s_{n_{|\mathcal{A}|}})=\tau_{n_{|\mathcal{A}|}}}\Lambda_F(\sum_{i=1}^{|\mathcal{A}|} s_i)\prod_{i=1}^{|\mathcal{A}|}G(s;\alpha_i)ds_i,\]
where, again, the contours are well chosen and $\mathcal{A}=\{n_1,\dots,n_{\mathcal{A}}\}$. Using this, one reduces the problem to slight generalizations of the techniques of Booker to deduce a contradication based on the location of a pole of $\Lambda_F(s)$. In fact, it is possible to deduce that the $L$-function of an $L$-data of degree $1<d<2$ can not have finite abscissa of convergence. This is analogous to the idea exploited in \cite{OTSOSCVII} that in the Selberg class $L$-functions are holomorphic in the right-plane $\Re(s)>1$. In our setting of exponential sums, holomorphy amounts to a notion of \textsl{cuspidality} and we detect poles via the constant term in a Fourier series.
\newline The following lemma is quickly deduced from the functional equation of $\Lambda_F(s)$ and is to be compared with \cite[section~3.2]{LFAD} and \cite[theorem~1.1]{OTSOSCVII}. 
\begin{lemma}
Let $F$ be an $L$-datum od degree $1<d<2$ and $L$-function $L_F(s)$. If $G(s;\alpha_i)$ is defined as above and $\mathcal{A}\subset\{1,\dots,N\}$, then
\[\sum_{\emptyset\neq\mathcal{A}\subset[N]}\frac{1}{(2\pi i)^{|\mathcal{A}|}}\int_{\Re(s_{n_1})=\tau_{n_1}}\cdots\int_{\Re(s_{n_{|\mathcal{A}|}})=\tau_{n_{|\mathcal{A}|}}}\Lambda_F(\sum_{i=1}^{|\mathcal{A}|}s_i)\prod_{i=1}^{|\mathcal{A}|}G(s;\alpha_i)ds_i\]
\[=O(\Im(z)^{-\varepsilon})+\frac{1}{(2\pi i)^{|\mathcal{A}|}}\int_{\Re(s_{n_1})=\tau_{n_1}}\cdots\int_{\Re(s_{n_{|\mathcal{A}|}})=\tau_{n_{|\mathcal{A}|}}}\sum_{n=1}^{\infty}\frac{\overline{a_F(n)}}{n^{1-\sum_{i=0}^Ns_i}}\cdot\]
\[\cdot\prod_{i=1}^{|\mathcal{A}|}A^{s/\alpha-\frac{1}{2}}\frac{\Gamma((s-\frac{1}{2})(d-\frac{1}{\alpha})+1-\frac{1}{2\alpha})}{\cos(\frac{\pi}{2}((s-\frac{1}{2})(\frac{2}{\alpha}-d)+\frac{1-\alpha+\mu\alpha}{\alpha}))}(-iz)^{s_i-1}ds_i,\]
for some $A\in\mathbb{R}_{>0}$
\end{lemma}
Though already somewhat elaborate, this lemma in fact admits a vast generalisation. The key is to understand the equation above as living in a family which can be described in terms of two operators acting on a class of functions as in \cite[Theorems~1.2,~1.3]{OTSOSCVII}. Applying a carefully chosen combination of these operators allows one deduce the existence of poles for $L_F(s)$ existing further and further to the right in the complex plane. 
\bibliography{RIMSrefs}

\begin{thebibliography}{10}

\bibitem{ITLP}
J.~Bernstein and S.~Gelbart.
\newblock {\em An Introduction to the {L}anglands Program}.
\newblock Birkh\"{a}user Basel, 2004.

\bibitem{BockleThesis}
G.~B\"{o}ckle.
\newblock Universal deformations of even {G}alois representations and relations
  to {M}aass waveforms.
\newblock PhD thesis, University of Illinois at Urbana-Champaign, 1995.

\bibitem{LFAD}
A.~Booker.
\newblock L-functions as distributions.
\newblock {\em Math. Ann}, 363:423--454, 2015.

\bibitem{WCTWP}
A.~Booker and M.~Krishnamurthy.
\newblock Weil's converse theorem with poles.
\newblock {\em International Mathematics Research Notices}, 127:1--12, 2013.

\bibitem{AFAR}
D.~Bump.
\newblock {\em Automorphic {F}orms and {R}epresentations}.
\newblock Cambridge Studies in Avanced Mathematics, 1998.

\bibitem{CTFGLn}
J.W. Cogdell and I.I. Piatetski-Shapiro.
\newblock Converse theorems for $ gl\_n$.
\newblock {\em Publications Math{\'e}matiques de l'IHES}, 79:157--214, 1994.

\bibitem{CTFGLnII}
J.W. Cogdell and I.I. Piatetski-Shapiro.
\newblock Converse theorems for $ gl\_n$ {II}.
\newblock {\em J. Reine Angew. Math.}, 507:165--188, 1999.

\bibitem{OTSOSCV}
J.~Kaczorowski and A.~Perelli.
\newblock On the structure of the {S}elberg class {V}. $1<d<5/3$.
\newblock {\em Invent. Math.}, 150:485--516, 2002.

\bibitem{OTSOSCVI}
J.~Kaczorowski and A.~Perelli.
\newblock On the structure of the {S}elberg class {VI}: non-linear twists.
\newblock {\em Acta. Arith.}, 116:315--341, 2005.

\bibitem{OTSOSCVII}
J.~Kaczorowski and A.~Perelli.
\newblock On the structure of the {S}elberg class {VII}. $1<d<2$.
\newblock {\em Ann. of Math}, 173:1397--1441, 2011.

\bibitem{COLDLD}
T.~Oliver.
\newblock Classification of low degree {L}-data.
\newblock in preparation, 2016.

\bibitem{ACOZOLF}
R.~Raghunathan.
\newblock A comparison of zeros of {L}-functions.
\newblock {\em Math. Res. Lett.}, 6:155--167, 1999.

\bibitem{CTDSP}
R.~Raghunathan.
\newblock A converse theorem for {D}irichlet series with poles, in: Cohomology
  of arithmetic groups, {L}-functions and automorphic forms, mumbai, 1998/1999.
\newblock {\em Tata Inst. Fund. Res. Stud. Math.}, 15:127--142, 2001.

\bibitem{OLFWPSMFE}
R.~Raghunathan.
\newblock On {L}-functions with poles satisfying {M}aas's functional equation.
\newblock {\em J. Number Theory}, 6:1255--1273, 2010.

\bibitem{ZODSWPC}
E.~Saias and A.~Weingartner.
\newblock Zeros of {D}irichlet series with periodic coefficients.
\newblock {\em Acta. Arith.}, 140:367--385, 2009.

\bibitem{OANCARAACODS}
A.~Selberg.
\newblock Old and new conjectures and results about a class of {D}irichlet
  series.
\newblock {\em Univ. Salerno}, Proceedings of the Amalfi Conference on Analytic
  Number Theory (Maiori, 1989):367--385, 1992.

\bibitem{SLFEDLTDNP}
A.~Weil.
\newblock Sur les ``formules explicites'' de la theorie des nombres premiers.
\newblock {\em Comm. S\'{e}m. Math. Univ. Lund [Medd. Lunds Univ. Mat. Sem.]},
  Tome Supplementaire:252--265, 1952.

\bibitem{UDBDRDF}
A.~Weil.
\newblock \"{U}ber die bestimmung dirichletscher reihen durch
  funktionalgleichungen.
\newblock {\em Mathematische Annalen}, 168:149--156, 1967.

\end{thebibliography}
\bibliographystyle{plain}
\end{document}